\documentclass[11pt,reqno]{amsart}
\allowdisplaybreaks[1]
\begin{document}

\theoremstyle{definition}
\newtheorem{defn}{Definition}
\newtheorem{lemma}{Lemma}
\newtheorem{prop}{Proposition}
\newtheorem{step}{Step}

\title[Bounds on ternary square-free words]{Improved 
bounds on the number of ternary square-free words}
\author{Uwe Grimm}
\address{Applied Mathematics Department, 
Faculty of Mathematics and Computing,
The Open University, Walton Hall, 
Milton Keynes MK7 6AA, U.K.}
\email{u.g.grimm@open.ac.uk}
\urladdr{http://mcs.open.ac.uk/ugg2}
\begin{abstract}
  Improved upper and lower bounds on the number of square-free ternary
  words are obtained. The upper bound is based on the enumeration of
  square-free ternary words up to length $110$. The lower bound is
  derived by constructing generalised Brinkhuis triples. The problem
  of finding such triples can essentially be reduced to a
  combinatorial problem, which can efficiently be treated by computer.
  In particular, it is shown that the number of square-free ternary
  words of length $n$ grows at least as $65^{n/40}$, replacing the
  previous best lower bound of $2^{n/17}$.
\end{abstract}
\maketitle

\section{Introduction}

A word $w$ is a string of letters from a certain alphabet $\Sigma$,
the number of letters of $w$ is called the length of the word. The set
of words of length $n$ is $\mathcal{L}(n)=\Sigma^n$, and the union
\begin{equation}
\mathcal{L}=\bigcup_{n\ge 0}\mathcal{L}(n)=\Sigma^{\mathbb{N}_0}
\end{equation}
is called the language of words in the alphabet $\Sigma$. This is a
monoid with concatenation of words as operation and the empty word
$\lambda$, which has zero length, as neutral element \cite{L}. For a
word $w$, we denote by $\bar{w}$ the corresponding reversed word,
i.e., the word obtained by reading $w$ from back to front. A
palindrome is a word $w$ that is symmetric, $w=\bar{w}$.

Square-free words [1--13]
%\cite{BEG,BEN,Bra,Bri,Cro1,Cro2,EZ,F,Kob,Lec,L,Ple,See} 
are words $w$ that do not contain a ``square'' $yy$ of a word $y$ as a
subword (factor). In other words, $w$ can only be written in the form
$xyyz$, with words $x$, $y$ and $z$, if $y=\lambda$ is the empty
word. In a two-letter alphabet $\{0,1\}$, the complete list of
square-free words is $\{\lambda,0,1,01,10,010,101\}$. However, in a
three-letter alphabet $\Sigma=\{0,1,2\}$, square-free words of
arbitrary length exist, and the number of square-free words of a given
length $n$ grows exponentially with $n$ \cite{Bri,Bra,EZ}.

We denote the set of square-free words of length $n$ in the alphabet
$\Sigma=\{0,1,2\}$ by $\mathcal{A}(n)\subset\mathcal{L}(n)$. The
language of ternary square-free words is 
\begin{equation}
\mathcal{A}=\bigcup_{n\ge 0}\mathcal{A}(n)\subset\Sigma^{\mathbb{N}_0}.
\end{equation}
We are interested in the number of square-free
words of length $n$
\begin{equation}
a(n)=|\mathcal{A}(n)|
\end{equation}
and in estimating the growth of $a(n)$ with the length $n$. 
For $n=0,1,2,3$, the sets of ternary square-free words are 
\begin{align}
\mathcal{A}(0) & = \{\lambda\},\\
\mathcal{A}(1) & = \{0,1,2\},\\
\mathcal{A}(2) & = \{01,02,10,12,20,21\},\\
\mathcal{A}(3) & = \{010,012,020,021,101,102,120,121,201,202,210,212\},
\end{align}
where $\lambda$ denotes the empty word. Hence $a(0)=1$, $a(1)=3$,
$a(2)=6$, $a(3)=12$, and so on, see \cite{BEG} where the values of
$a(n)$ for $n\le 90$ are tabulated. In \cite{SP}, the sequence is
listed as A006156 (formerly M2550).

\section{Upper bounds obtained by enumeration}

Obviously, a word $w$ of length $m+n$, obtained by concatenation of
words $w_1$ of length $m$ and $w_2$ of length $n$, can only be
square-free if $w_1$ and $w_2$ are square-free. This necessary, but not
sufficient, condition implies the inequality
\begin{equation}
a(m+n)\le a(m)a(n)
\label{ub1}
\end{equation}
for all $m,n\ge 0$. By standard arguments, see also \cite{BEG}, this
guarantees the existence of the limit
\begin{equation}
s:=\lim_{n\rightarrow\infty} a(n)^{\frac{1}{n}},
\label{s}
\end{equation}
the growth rate or ``connective constant'' of ternary square-free
words \cite{F}. The precise value of $s$ is not known, but lower
\cite{Bri,Bra,EZ} and upper bounds \cite{BEG} have been established.
It is the purpose of this paper to improve both the lower and the
upper bounds.

It is relatively easy to derive reasonable upper bounds from the
inequality (\ref{ub1}). In fact \cite{BEG}, one can slightly improve
on (\ref{ub1}) by considering two words $w_1$ and $w_2$ of length
$m\ge 2$ and $n\ge 2$, such that the last two letters of $w_1$ are
equal to the first two letters of $w_2$, and we join them to a word
$w$ of length $m+n-2$ by having the two words overlap on these two
letters. This yields
\begin{equation}
a(m+n-2)\le\frac{1}{6}a(m)a(n),
\label{ub2}
\end{equation}
for all $m,n\ge 2$, because there are precisely $a(n)/6$ square-free
letters of length $n\ge 2$ that start with the last two letters of
$w_1$. Taking $n$ fixed, one obtains
\begin{equation}
s^{n-2}=\lim_{m\rightarrow\infty}\frac{a(m+n-2)}{a(m)}
\le\frac{a(n)}{6}
\label{ub3}
\end{equation}
and hence the upper bound
\begin{equation}
s\le\left(\frac{a(n)}{6}\right)^{\frac{1}{n-2}}
\label{ub4}
\end{equation}
for any $n\ge 3$. This bound can be systematically improved by
calculating $a(n)$ for as large values of $n$ as possible. The bound
given in \cite{BEG}, from $a(90)=258\,615\,015\,792$, is
\begin{equation}
s\le{43\,102\,502\,632}^{\frac{1}{88}}=1.320\,829\,\ldots
\end{equation}
The results given in table~\ref{sfw} extend the previously known
values of $a(n)$ \cite{BEG} to lengths $n\le 110$. They were obtained
by a simple algorithm, extending square-free words letter by letter
and checking that the new letter does not lead to the formation of any
square. The value $a(110)$ yields an improved upper bound of
\begin{equation}
s\le {8\,416\,550\,317\,984}^{\frac{1}{108}}=1.317\,277\,\ldots
\end{equation}

\begin{table}
\caption{The number of ternary square-free words $a(n)$ of length $n$
for $91\le n\le 110$.\label{sfw}}
\begin{center}
\begin{tabular}{@{}r@{\extracolsep{2ex}}r@{\extracolsep{6ex}}r@{\extracolsep{2ex}}r@{}}
\hline
\multicolumn{1}{c}{$n$} &
\multicolumn{1}{c}{$a(n)$} & 
\multicolumn{1}{c}{$n$} &
\multicolumn{1}{c}{$a(n)$}\\
\hline
  $91$ &     $336\,655\,224\,582$ & 
 $101$ &  $4\,704\,369\,434\,772$ \\
  $92$ &     $438\,245\,025\,942$ & 
 $102$ &  $6\,123\,969\,129\,810$ \\
  $93$ &     $570\,491\,023\,872$ &
 $103$ &  $7\,971\,950\,000\,520$ \\
  $94$ &     $742\,643\,501\,460$ &
 $104$ & $10\,377\,579\,748\,374$ \\
  $95$ &     $966\,745\,068\,408$ &
 $105$ & $13\,509\,138\,183\,162$ \\
  $96$ &  $1\,258\,471\,821\,174$ &
 $106$ & $17\,585\,681\,474\,148$ \\
  $97$ &  $1\,638\,231\,187\,596$ &
 $107$ & $22\,892\,370\,891\,330$ \\
  $98$ &  $2\,132\,586\,986\,466$ &
 $108$ & $29\,800\,413\,809\,730$ \\
  $99$ &  $2\,776\,120\,525\,176$ &
 $109$ & $38\,793\,041\,799\,498$ \\
 $100$ &  $3\,613\,847\,436\,684$ &
 $110$ & $50\,499\,301\,907\,904$ \\
\hline
\end{tabular}
\end{center}
\end{table}

\section{Brinkhuis triples and lower bounds}

While the upper bound is already relatively close to the actual value
of $s$, which was estimated in reference \cite{BEG} to be about
$1.301\,76$ on the basis of the first $90$ values, it is much more
difficult to obtain any reasonable lower bound for $s$.  In order to
derive a lower bound, one has to show that $a(n)$ grows exponentially
in $n$ with optimal growth bound. This can be achieved by
demonstrating that each square-free word of length $n$ gives rise to
sufficiently many different square-free words of some length
$m>n$. This was first done by Brinkhuis \cite{Bri}, by constructing
what is now known as a Brinkhuis triple or a Brinkhuis triple pair.

\begin{defn}
  An $n$-Brinkhuis triple pair is a set
  $\mathcal{B}=\{\mathcal{B}^{(0)},\mathcal{B}^{(1)},\mathcal{B}^{(2)}\}$
  of three pairs
  $\mathcal{B}^{(i)}=\{U^{(i)},V^{(i)}\}\subset\mathcal{A}(n)$,
  $i\in\{0,1,2\}$, of pairwise different square-free words such that
  the set of $96$ words of length $3n$
\[
  \bigcup_{w_1w_2w_3\in\mathcal{A}(3)}
  \left\{W_1W_2W_3\mid W_j\in\mathcal{B}^{(w_j)},
  j=1,2,3\right\}\subset\mathcal{A}(3n).
\]
\label{bt}
\end{defn}

In other words, it is required that all $3n$-letter images of the
twelve elements of $\mathcal{A}(3)$ under any combination of the eight
maps
\begin{equation}
\varrho_{x,y,z}:
\begin{cases}
0\rightarrow x\in\mathcal{B}^{(0)}\\
1\rightarrow y\in\mathcal{B}^{(1)}\\
2\rightarrow z\in\mathcal{B}^{(2)}
\end{cases} 
\end{equation}
are square-free. This property is sufficient to ensure that images of
any square-free word in the alphabet $\Sigma$ under any combination
of the eight maps to each of its letters is again square-free.  This
can be shown as follows.

Consider the six-letter alphabet $\tilde{\Sigma}= \{0,0',1,1',2,2'\}$
and a language $\tilde{\mathcal{A}}$ consisting of all words of
$\mathcal{A}$ with an arbitrary number of letters replaced by their
primed companions. In other words,
\begin{equation}
\tilde{\mathcal{A}}=\bigcup_{n\ge 0}\tilde{\mathcal{A}}(n),\qquad
\tilde{\mathcal{A}}(n) = \left\{w\in\tilde{\Sigma}^{n}\mid
\pi(w)\in\mathcal{A}(n)\right\}
\end{equation}
where $\pi$ is the map 
\begin{equation}
\pi:\tilde{\Sigma}\rightarrow\Sigma,\quad
\pi(0)\!=\!\pi(0')\!=\!0,\; 
\pi(1)\!=\!\pi(1')\!=\!1,\; 
\pi(2)\!=\!\pi(2')\!=\!2, 
\end{equation}
that projects back to the three-letter alphabet $\Sigma$. 
The map
\begin{equation}
\varrho:
\begin{cases}
0\rightarrow  U^{(0)},\; 0'\rightarrow V^{(0)}\\
1\rightarrow  U^{(1)},\; 1'\rightarrow V^{(1)}\\
2\rightarrow  U^{(2)},\; 2'\rightarrow V^{(2)}
\end{cases}
\end{equation}
is a uniformly growing morphism from the language
$\tilde{\mathcal{A}}$ into the language $\mathcal{L}$. By the
condition (\ref{bt}), this morphism is square-free on all three-letter
words in $\tilde{\mathcal{A}}$, i.e., the images of elements in
$\tilde{\mathcal{A}}(3)$ are square-free.  As $\varrho$ is a uniformly
growing morphisms, being square-free on $\tilde{\mathcal{A}}(3)$
implies, as proven in \cite{Cro1} and \cite{Bra}, that $\varrho$ is a
square-free morphism, i.e., it maps square-free words in
$\tilde{\mathcal{A}}$ onto square-free words in $\mathcal{L}$, thus
onto words in $\mathcal{A}$.

\begin{lemma}
  The existence of an $n$-Brinkhuis triple pair implies the lower
  bound $s\ge 2^{1/(n-1)}$.
\label{lem}
\end{lemma}
\begin{proof}
  The existence of an $n$-Brinkhuis triple pair implies the inequality
  \begin{equation} 
  a(mn)\ge 2^m a(m)
  \end{equation}
  for any $m>0$, because each square-free word of length $m$ yields
  $2^m$ different square-free words of length $mn$. This means
  \begin{equation}
  \left(\frac{a(mn)}{a(m)}\right)^{\frac{1}{m}}\ge 2,
  \end{equation}
  for any $m>0$, and hence
  \begin{equation}
  s^{n-1}=\lim_{m\rightarrow\infty}
  \left(\frac{a(mn)}{a(m)}\right)^{\frac{1}{m}}\ge 2,
  \end{equation}
  establishing the lower bound.
\end{proof}

The first lower bound was derived by Brinkhuis \cite{Bri}, who showed
that $s\ge 2^{1/24}$ by constructing a $25$-Brinkhuis triple pair
consisting entirely of palindromic words. In that case, the conditions
on the square-freeness of the images of three-letter words can be
simplified to the square-freeness of the images of two-letter words
and certain conditions on the ``heads'' and the ``tails'' of the
words, which are easier to check explicitly. Brandenburg \cite{Bra}
produced a $22$-Brinkhuis triple pair, which proves a lower bound of
$s\ge 2^{1/21}$. For a long time, this was the best lower bound
available, until, quite recently, Ekhad and Zeilberger \cite{EZ} came
up with a $18$-Brinkhuis triple pair equivalent to
\begin{align}
U^{(0)} & = 012021020102120210 &
V^{(0)} & = 012021201020120210 = \bar{U}^{(0)} \notag\\[-1ex]
U^{(1)} & = 120102101210201021 &
V^{(1)} & = 120102012101201021 = \bar{U}^{(1)} \notag\\[-1ex]
U^{(2)} & = 201210212021012102 & 
V^{(2)} & = 201210120212012102 = \bar{U}^{(2)},
\label{EZtrip}
\end{align}
thus establishing the bound $s\ge 2^{1/17}$.  We note that the simpler
definition for a Brinkhuis triple pair in \cite{EZ}, which is akin to
Brinkhuis' original approach, is in fact incomplete, as it does not rule out a
square that overlaps three adjacent words if the words are not palindromic.
Nevertheless, the Brinkhuis triple (\ref{EZtrip}) given in \cite{EZ} is
correct, and so is the lower bound $s\ge 2^{1/17}=1.041\,616\,\ldots$ derived
from it. In fact, it has been claimed (see \cite{Z}) that this is the optimal
bound that can be obtained in this way, and this is indeed the case, see the
discussion below.

It is interesting to note that, although this minimal-length Brinkhuis
triple pair does not consist of palindromes, it is nevertheless
invariant under reversion of words, as $V^{(i)}=\bar{U}^{(i)}$. In
addition, it also shares the property with Brinkhuis' orginial triple
that the words $U^{(1)}$, $V^{(1)}$ and $U^{(2)}$, $V^{(2)}$ which
replace the letters $1$ and $2$, respectively, are obtained from
$U^{(0)}$, $V^{(0)}$ by a global permutation $\tau$ of the three
letters
\begin{equation}
\tau: 
\begin{cases}
0 \rightarrow 1\\[-0.5ex]
1 \rightarrow 2\\[-0.5ex]
2 \rightarrow 0
\end{cases}
\label{tau}
\end{equation}
i.e.,
\begin{equation}
U^{(2)}=\tau(U^{(1)})=\tau^{2}(U^{(0)}),\quad
V^{(2)}=\tau(V^{(1)})=\tau^{2}(V^{(0)}).
\end{equation}
Clearly, given any Brinkhuis triple pair, the sets of words obtained
by reversion or by applying any permutation of the letters are again
Brinkhuis triple pairs, so it may not be too surprising that a
Brinkhuis triple pair of minimal length turns out to be invariant
under these two operations.

\section{Generalised Brinkhuis triples}

As we cannot improve on the lower bound by constructing a shorter
Brinkhuis triple, we proceed by generalising the notion. The idea is
to allow for more than two words that replace each letter. This leads
to the following general definition.
\begin{defn}
  An $n$-Brinkhuis $(k_0,k_1,k_2)$-triple is a set of $k_0+k_1+k_2$
  square-free words
  $\mathcal{B}=\{\mathcal{B}^{(0)},\mathcal{B}^{(1)},\mathcal{B}^{(2)}\}$,
  $\mathcal{B}^{(i)}=\{w^{(i)}_{j}\in S(n)\mid 1\le j\le k_i\}$, $k_i\ge 1$,
  such that,
  for any square-free word $ii'i''$ of length $3$ and any $1\le j\le k_i$,
  $1\le j'\le k_{i'}$, $1\le j''\le k_{i''}$, the composed word
  $w^{(i)}_{j}w^{(i')}_{j'}w^{(i'')}_{j''}$ of length $3n$ is square-free.
  \label{gentrip}
\end{defn}
Note that the definition reduces to definition~\ref{bt} in the case
$k_0=k_1=k_2=2$ of an ``ordinary'' Brinkhuis triple pair. {}From the
set of square-free words of length $3$, we deduce that the number of
composed words that enter is
$6k_{0}k_{1}k_{2}+k_{0}^2(k_{1}+k_{2})+k_{1}^2(k_{0}+k_{2})+
k_{2}^2(k_{0}+k_{1})$.

\begin{lemma}
  The existence of an $n$-Brinkhuis $(k_0,k_1,k_2)$-triple implies the lower
  bound $s\ge k^{1/(n-1)}$, where $k=\min(k_0,k_1,k_2)$.
\end{lemma}
\begin{proof}
  The proof proceeds as in lemma~\ref{lem} above, with $2$ replaced by
  $k=\min(k_0,k_1,k_2)$.
\end{proof}

As far as the lower bound is concerned, we do not gain anything by
considering triples where the number of words $k_0$, $k_1$ and $k_2$
differ from each other. Nevertheless, the generality of
definition~\ref{gentrip} shall be of use below. In order to derive
improved lower bounds, we shall in fact concentrate on a more
restricted class of triples.

\begin{defn}
  A special $n$-Brinkhuis $k$-triple is an $n$-Brinkhuis $(k,k,k)$-triple
  $\mathcal{B}=\{\mathcal{B}^{(0)},\mathcal{B}^{(1)},\mathcal{B}^{(2)}\}$
  such that
  $\mathcal{B}^{(2)}=\tau(\mathcal{B}^{(1)})=\tau^{2}(\mathcal{B}^{(0)})$
  and $w\in\mathcal{B}^{(0)}$ implies $\bar{w}\in\mathcal{B}^{(0)}$, where
  $\tau$ is the permutation of letters defined in equation (\ref{tau}).
\end{defn}

The first condition means that all words in $\mathcal{B}^{(1)}$ and
$\mathcal{B}^{(2)}$ can be obtained from the words in
$\mathcal{B}^{(0)}$ by the global permutation $\tau$. The second
condition implies that the words in $\mathcal{B}^{(0)}$, and hence
also in $\mathcal{B}^{(1)}$ and $\mathcal{B}^{(2)}$, are either
palindromes, i.e., $w=\bar{w}$, or occur as pairs $(w,\bar{w})$.  This
means that a special $n$-Brinkhuis $k$-triple is characterised by the
set of palindromes $w=\bar{w}\in\mathcal{B}^{(0)}$ and by one member
of each pairs of non-palindromic words
$(w,\bar{w})\in\mathcal{B}^{(0)}$. If there are $k_{\rm p}$
palindromes and $k_{\rm n}$ pairs in $\mathcal{B}^{(0)}$, then these
generate a special Brinkhuis $k$-triple with $k=k_{\rm p}+2k_{\rm
n}$. We shall call $K=(k_{\rm p},k_{\rm n})$ the signature of the
special Brinkhuis $k$-triple, and denote a set of $k_{\rm p}+k_{\rm
n}$ generating words by $\mathcal{G}$.

In order to obtain the best lower bound possible, we are looking for
optimal choices of the length $n$ and the number of words $k$. There
are two possibilities, we may look for the largest $k$ for given
length $n$, or for the smallest length $n$ for a given number
$k$. This is made precise by the following definitions.

\begin{defn}
  An optimal special $n$-Brinkhuis triple is a special $n$-Brinkhuis
  $k$-triple such that any special $n$-Brinkhuis $l$-triple has $l\le k$. 
\end{defn}

\begin{defn}
  A minimal-length special Brinkhuis $k$-triple is a special $n$-Brinkhuis
  $k$-triple such that any special $m$-Brinkhuis $k$-triple has $m\ge n$. 
\end{defn}

If $\mathcal{B}$ is a special $n$-Brinkhuis $k$-triple, so is its
image $\sigma(\mathcal{B})$ under any permutation $\sigma\in S_{3}$ of
the three letters. Therefore, without loss of generality, we may
assume that the first word $w^{(0)}_{1}\in\mathcal{B}$ starts with the
letters $01$. This has the following consequences on the other words
of the triple.

\begin{lemma}
  Consider a special $n$-Brinkhuis $k$-triple $\mathcal{B}$, with
  $n>1$, such that the word $w^{(0)}_{1}\in\mathcal{B}^{(0)}$ starts
  with the letters $01$.  Then $n\ge 7$, and all words in
  $\mathcal{B}^{(0)}$ start with the three letters $012$ and end on
  $210$.  \label{ht1}
\end{lemma}
\begin{proof}
  As $w^{(1)}_{1}=\tau(w^{(0)}_{1})$ and
  $w^{(2)}_{1}=\tau^{2}(w^{(0)}_{1})$, the words $w^{(1)}_{1}$ and
  $w^{(2)}_{1}$ start with letters $12$ and $20$, respectively.  If
  $n=2$, then $w^{(0)}_{1}w^{(1)}_{1}=0112$ contains the square $11$,
  so $n\ge 3$.  Square-freeness of the composed words
  $w^{(0)}_{j}w^{(1)}_{1}$ and $w^{(0)}_{j}w^{(2)}_{1}$, $1\le j\le
  k$, implies that the words $w^{(0)}_{j}$ have to end on $210$,
  because $w=210$ is the only word in $\mathcal{A}(3)$ such that $w12$
  and $w20$ are both square-free.  This in turn implies that all words
  in $\mathcal{B}^{(1)}$ and $\mathcal{B}^{(2)}$ end on $021$ and
  $102$, respectively. Now, square-freeness of the composed words
  $w^{(1)}_{j}w^{(0)}_{j'}$ and $w^{(2)}_{j}w^{(0)}_{j'}$ implies that
  the first three letter of $w^{(0)}_{j}$, for any $1\le j\le k$, have
  to be $w=012$, because this is the only word $\mathcal{A}(3)$ in
  such that $021w$ and $102w$ are both square-free. For $n=3$ and
  $n=4$, no such words exist, and the only possibility for $n=6$ would
  be $012210$ which is not square-free.  For $n=5$, the square-free
  word $01210$ starts with $012$ and ends on $210$, but
  $w^{(0)}_{1}w^{(2)}_{1}w^{(0)}_{1}=012102010201210$ contains the
  square of $0201$.
\end{proof}

One can even say more about the ``heads'' and ``tails'' of the words
in a special Brinkhuis triple. There are two possible choices for the
forth letter of $w^{(0)}_{1}$, and both possibilities fix further 
letters and cannot appear within the same special Brinkhuis triple.
Therefore, we can distinguish two different types of special Brinkhuis
triples.
\begin{prop}
  Consider a special $n$-Brinkhuis $k$-triple $\mathcal{B}$, with
  $n>1$, such that the word $w^{(0)}_{1}\in\mathcal{B}^{(0)}$ starts
  with the letters $01$.  Then $n\ge 13$ and either all words in
  $\mathcal{B}^{(0)}$ are of the form $012021\ldots 120210$, or all
  words are of the form $012102\ldots 201210$.
\label{ht2}
\end{prop}
\begin{proof}
  {}From lemma~\ref{ht1}, we know that $n\ge 7$ and $w^{(0)}_{1}$
  starts with $012$ and ends on $210$. There are now two choices for
  the forth letter.  Let us consider the case that $w^{(0)}_{1}$
  starts with $0120$. Then $w^{(1)}_{1}$ starts with $1201$. Now, from
  lemma~\ref{ht1}, $w^{(2)}_{1}$ ends on $102$, and square-freeness of
  $w^{(2)}_{1}w^{(0)}_{j}$ implies that $w^{(0)}_{j}$ starts with
  $01202$, and hence with $012021$. Now $w^{(1)}_{1}$ starts with
  $120102$ and $w^{(2)}_{1}$ with $201210$. From square-freeness of
  $w^{(0)}_{j}w^{(1)}_{1}$ and $w^{(0)}_{j}w^{(2)}_{1}$, we can rule
  out $w^{(0)}_{j}$ ends on $1210$, because both possible extension
  $101210$ and $201210$ result in squares.  Hence $w^{(0)}_{j}$ ends
  on $0210$ and, from square-freeness of $w^{(0)}_{j}w^{(w)}_{1}$,
  it has to end on $20210$, and thus on $120210$.
  
  Consider now the second possibility, i.e., $w^{(0)}_{1}$ starts with
  $0121$.  Necessarily, it then starts with $01210$. As $w^{(1)}_{1}$
  ends on $021$, square-freeness of $w^{(1)}_{1}w^{(0)}_{j}$ implies
  that $w^{(0)}_{j}$ starts with $012102$. Then $w^{(1)}_{1}$ starts
  with $120210$ and $w^{(2)}_{1}$ with $201021$. Square-freeness of
  $w^{(0)}_{j}w^{(1)}_{1}$ and $w^{(0)}_{j}w^{(2)}_{1}$ rules out an
  ending $0210$ for $w^{(0)}_{j}$, as the only possible extension
  $20120$ and $120210$ both result in squares.  Hence $w^{(0)}_{j}$
  ends on $01210$, and, from square-freeness of
  $w^{(0)}_{j}w^{(1)}_{1}$, actually has to end on $201210$.
  
  Now, in both cases it is obviously impossible to find square-free
  words of length $n=8,9,10,12$ that satisfy these conditions. For the
  first case, the one choice left for $n=11$ is $01202120210$, which
  contains the square of $1202$. In the second case, the only word for
  $n=11$ that satisfies the conditions is $01210201210$. In this case,
  $w^{(0)}_{1}w^{(1)}_{1}=0121020121012021012021$ contains the square
  of $210120$.
\end{proof}

The proofs of lemmas~\ref{ht1} and \ref{ht2} are very explicit, but
you may simplify the argument by realising that the conditions at both
ends are essentially equivalent, as they follow from reversing the
order of letters in combined words. The results restrict the number of
words that have to be taken into account when looking for a special
$n$-Brinkhuis $k$-triple. In what follows, we can restrict ourselves
to the case $n\ge 13$. We denote the set of such square-free words by
\begin{align}
\mathcal{A}_{1}(n)&=\{w\in\mathcal{A}(n)\mid w=012021\ldots120210\}
\subset\mathcal{A}(n),\label{A1}\\
\mathcal{A}_{2}(n)&=\{w\in\mathcal{A}(n)\mid w=012102\ldots201210\}
\subset\mathcal{A}(n)\label{A2},
\end{align}
and the number of such words by
\begin{align}
a_1(n)&:=|\mathcal{A}_{1}(n)|,\\
a_2(n)&:=|\mathcal{A}_{2}(n)|.
\end{align}
We denote the number of palindromes by
\begin{align}
a_{1{\rm p}}(n)&:=|\{w\in\mathcal{A}_{1}(n)\mid w=\bar{w}\}|,\\
a_{2{\rm p}}(n)&:=|\{w\in\mathcal{A}_{2}(n)\mid w=\bar{w}\}|,
\end{align}
and the number of non-palindromic pairs by
\begin{align}
a_{1{\rm n}}(n)&:=\frac{1}{2}(a_1(n)-a_{1{\rm p}}(n)),\\
a_{2{\rm n}}(n)&:=\frac{1}{2}(a_2(n)-a_{2{\rm p}}(n)).
\end{align}
Clearly, there are no palindromic square-free words of even length,
and thus $a_{1{\rm p}}(2n)=a_{2{\rm p}}(2n)=0$,  
$a_{1{\rm n}}(2n)=a_{1}(2n)/2$ and $a_{2{\rm n}}(2n)=a_{2}(2n)/2$.

Now, for a word $w\in\mathcal{A}_{1}(n)$ or $w\in\mathcal{A}_{2}(n)$
to be a member of a special $n$-Brinkhuis triple, it must at least
generate a triple by itself. This motivates the following definition.

\begin{defn}
  A square-free palindrome $w=\bar{w}\in\mathcal{A}(n)$ is called
  admissible if $w$ generates a special $n$-Brinkhuis $1$-triple. A
  non-palindromic square-free word $w$ of length $n$ is admissible if
  $w$ generates a special $n$-Brinkhuis $2$-triple.
\end{defn}

The hunt for optimal special $n$-Brinkhuis triples now proceeds in
three steps. 

\begin{step}
  The first step consists of selecting all admissible words in
  $\mathcal{A}_1(n)$ and $\mathcal{A}_2(n)$. Let us denote the number
  of admissible palindromes in $\mathcal{A}_1(n)$ by $b_{1{\rm p}}(n)$
  and the number of admissible non-palindromes by $2b_{1{\rm n}}(n)$,
  such that $b_{1{\rm n}}$ is the number of admissible pairs
  $(w,\bar{w})$ of non-palindromic words in
  $\mathcal{A}_1(n)$. Analogously, we define $b_{2{\rm p}}(n)$ and
  $b_{2{\rm n}}(n)$ for admissible words in $\mathcal{A}_2(n)$.
\end{step}
\begin{step}
  The second step consists of finding all triples of admissible words
  that generate a special $n$-Brinkhuis triple. Depending on the
  number of palindromes $k_{\rm p}$ in that triple, which can be
  $k_{\rm p}=0,1,2,3$, these are special Brinkhuis $k$-triples with
  $k=6,5,4,3$, respectively. We denote the number of such admissible
  triples by $t_{1}(n)$ and $t_{2}(n)$. Here, we need to check the
  conditions of definition~\ref{gentrip} for each triple. Using the
  structure of the special Brinkhuis triple, the number of words that
  have to be checked is substantially reduced from $12k^3$ to
  $k(2k^2+k_{\rm p})$.
\end{step}
\begin{step}
  The third and final step is purely combinatorial in nature, and does
  not involve any explicit checking of square-freeness of composed
  words. The reason is the following. A set $\mathcal{G}$,
  $|\mathcal{G}|\ge 3$, of words in $\mathcal{A}_1(n)$ or
  $\mathcal{A}_2(n)$, generates a special $n$-Brinkhuis triple if and
  only if all three-elemental subsets of $\mathcal{G}$ generate
  special $n$-Brinkhuis triples. This is obvious, because the
  conditions of definition~\ref{gentrip} on three-letter words never
  involve more than three words simultaneously, so checking the
  condition for all subsets of three generating words is necessary and
  sufficient. Thus, the task is to find the largest sets of generating
  words such that all three-elemental subsets are contained in our
  list of admissible triples. In order to obtain an optimal special
  $n$-Brinkhuis triple, one has to take into account that $k=k_{\rm
  p}+2k_{\rm n}$, so solutions with maximum number of generators are
  not necessarily optimal.
\end{step}

Even though this step is purely combinatorial and no further
operations on the words are required, it is by far the most expensive
part of the algorithm as the length $n$ increases. Therefore, this is
the part that limits the maximum length $n$ that we can consider.
Using a computer, we found the optimal Brinkhuis triples for $n\le
41$.  The results for generating words from $\mathcal{A}_{1}(n)$ are
given in table~\ref{trip1}, those for generating words taken from
$\mathcal{A}_{2}(n)$ are displayed in table~\ref{trip2}. We included
partial results for $42\le n\le 45$, in order to show how the number
of admissible words grows for larger $n$. Even though we do not know
the optimal $n$-Brinkhuis triples for these cases, it has to be
expected that the value of $k$ that can be achieved continues to grow,
and it is certainly true for $n=42$ where $k_{\rm opt}\ge 72$.

% Case 1: words of type $012021...120210$
% Columns:
% 
% length of words
% number of square-free words
% number of square-free words of form 012021...120210
% number of such palindromes
% number of pairs of non-palindromic such words
% number of eligible palindromes
% number of eligible non-palindromic pairs
% number of eligible triples of generators
% signature of optimal Brinkhuis triple
% optimal value of k 

\begin{table}
\caption{Results of the algorithm to find optimal special $n$-Brinkhuis 
triples with generating words in $\mathcal{A}_{1}(n)$.\label{trip1}}
\begin{tabular}{rrrrrrrrr@{,}lr}
\hline
$n$ &
$a$ &
$a_{1}$ &
$a_{1{\rm p}}$ &
$a_{1{\rm n}}$ &
$b_{1{\rm p}}$ &
$b_{1{\rm n}}$ &
$t_{1}$ &
\multicolumn{2}{c}{sign.} &
$k_{\rm opt}$\\
\hline
%n&          a&    a1&ap& an&bp& bn&       t& sign. & kopt\\
13&        342&     0& 0&  0& 0&  0&       0& (0&0) & 0\\
14&        456&     0& 0&  0& 0&  0&       0& (0&0) & 0\\
15&        618&     1& 1&  0& 0&  0&       0& (0&0) & 0\\
16&        798&     0& 0&  0& 0&  0&       0& (0&0) & 0\\
17&     1\,044&     1& 1&  0& 1&  0&       0& (1&0) & 1\\
18&     1\,392&     4& 0&  2& 0&  1&       0& (0&1) & 2\\
19&     1\,830&     5& 1&  2& 1&  0&       0& (0&0) & 0\\ 
20&     2\,388&     4& 0&  2& 0&  0&       0& (0&0) & 0\\
21&     3\,180&     1& 1&  0& 0&  0&       0& (0&0) & 0\\
22&     4\,146&     2& 0&  1& 0&  0&       0& (0&0) & 0\\
23&     5\,418&     3& 1&  1& 0&  1&       0& (0&1) & 2\\
24&     7\,032&     4& 0&  2& 0&  1&       0& (0&1) & 2\\
25&     9\,198&    13& 3&  5& 2&  1&       1& (2&1) & 4\\
26&    11\,892&    16& 0&  8& 0&  1&       0& (0&1) & 2\\
27&    15\,486&    18& 2&  8& 2&  0&       0& (1&0) & 1\\
28&    20\,220&    10& 0&  5& 0&  1&       0& (0&1) & 2\\
29&    26\,424&    27& 3& 12& 2&  3&       4& (2&2) & 6\\
30&    34\,422&    52& 0& 26& 0&  4&       0& (0&2) & 4\\
31&    44\,862&    64& 4& 30& 2&  7&       8& (1&3) & 7\\
32&    58\,446&    64& 0& 32& 0&  6&       5& (0&4) & 8\\
33&    76\,122&    60& 6& 27& 3&  7&      30& (0&6) &12\\
34&    99\,276&    70& 0& 35& 0&  7&      13& (0&4) & 8\\
35&   129\,516&   109& 9& 50& 4& 13&     328& (2&8) &18\\
36&   168\,546&   174& 0& 87& 0& 27&  1\,304& (0&15)&30\\
37&   219\,516&   291& 9&141& 6& 27&  2\,533& (3&14)&31\\
38&   285\,750&   376& 0&188& 0& 30&     973& (0&14)&28\\
39&   372\,204&   386&12&187& 3& 35&  2\,478& (2&15)&32\\
40&   484\,446&   428& 0&214& 0& 55& 10\,767& (0&24)&48\\
41&   630\,666&   593&15&289& 4& 76& 28\,971& (3&31)&65\\
42&   821\,154&   926& 0&463& 0&114& 74\,080& \multicolumn{2}{c}{?} & ?\\
43&1\,069\,512&1\,273&23&625&12&156&229\,180& \multicolumn{2}{c}{?} & ?\\
44&1\,392\,270&1\,518& 0&759& 0&170&235\,539& \multicolumn{2}{c}{?} & ?\\
45&1\,812\,876&1\,788&26&881&17&191&510\,345& \multicolumn{2}{c}{?} & ?\\
\hline
\end{tabular}
\end{table}
\clearpage

% Case 2: words of type $012102...201210$
% Columns:
% 
% length of words
% number of square-free words
% number of square-free words of form 012102...201210
% number of such palindromes
% number of pairs of non-palindromic such words
% number of eligible palindromes
% number of eligible non-palindromic pairs
% number of triples of generators

\begin{table}
\caption{Results of the algorithm to find optimal special $n$-Brinkhuis 
triples with generating words in $\mathcal{A}_{2}(n)$.\label{trip2}}
\begin{tabular}{rrrrrrrrr@{,}lr}
\hline
$n$ &
$a$ &
$a_{2}$ &
$a_{2{\rm p}}$ &
$a_{2{\rm n}}$ &
$b_{2{\rm p}}$ &
$b_{2{\rm n}}$ &
$t_{2}$ &
\multicolumn{2}{c}{sign.} &
$k_{\rm opt}$\\
\hline
%n&          a&  a2&ap& an&bp& bn&       t& sign. & kopt\\
13&        342&   1& 1&  0& 1&  0&       0& (1&0) & 1\\
14&        456&   0& 0&  0& 0&  0&       0& (0&0) & 0\\
15&        618&   0& 0&  0& 0&  0&       0& (0&0) & 0\\
16&        798&   0& 0&  0& 0&  0&       0& (0&0) & 0\\
17&     1\,044&   2& 0&  1& 0&  0&       0& (0&0) & 0\\
18&     1\,392&   2& 0&  1& 0&  0&       0& (0&0) & 0\\
19&     1\,830&   1& 1&  0& 0&  0&       0& (0&0) & 0\\
20&     2\,388&   0& 0&  0& 0&  0&       0& (0&0) & 0\\
21&     3\,180&   1& 1&  0& 0&  0&       0& (0&0) & 0\\
22&     4\,146&   6& 0&  3& 0&  0&       0& (0&0) & 0\\
23&     5\,418&   6& 2&  2& 2&  1&       0& (1&1) & 3\\
24&     7\,032&  10& 0&  5& 0&  2&       0& (0&1) & 2\\
25&     9\,198&  11& 1&  5& 1&  2&       1& (1&2) & 5\\
26&    11\,892&   8& 0&  4& 0&  1&       0& (0&1) & 2\\
27&    15\,486&   8& 2&  3& 1&  1&       0& (1&1) & 3\\
28&    20\,220&  10& 0&  5& 0&  3&       0& (0&2) & 4\\
29&    26\,424&  30& 4& 13& 1&  3&       2& (0&3) & 6\\
30&    34\,422&  40& 0& 20& 0&  6&       5& (0&4) & 8\\
31&    44\,862&  37& 5& 16& 2&  3&       2& (1&2) & 5\\
32&    58\,446&  32& 0& 16& 0&  4&       0& (0&2) & 4\\
33&    76\,122&  49& 5& 22& 2&  3&       7& (1&3) & 7\\
34&    99\,276&  76& 0& 38& 0& 10&      39& (0&5) &10\\ 
35&   129\,516& 142& 6& 68& 3& 20&     483& (2&7) &16\\
36&   168\,546& 188& 0& 94& 0& 29&  1\,602& (0&16)&32\\
37&   219\,516& 205& 9& 98& 3& 32&  2\,707& (1&13)&27\\
38&   285\,750& 198& 0& 99& 0& 27&  1\,112& (0&11)&22\\
39&   372\,204& 231&13&109& 6& 36&  5\,117& (2&14)&30\\
40&   484\,446& 396& 0&198& 0& 56& 12\,002& (0&19)&38\\
41&   630\,666& 615&15&300& 8& 81& 54\,340& (1&29)&59\\ 
42&   821\,154& 820& 0&410& 0&120&123\,610& \multicolumn{2}{c}{?} & ?\\
43&1\,069\,512& 969&15&477&10&158&332\,054& \multicolumn{2}{c}{?} & ?\\
44&1\,392\,270&1070& 0&535& 0&166&362\,560& \multicolumn{2}{c}{?} & ?\\
45&1\,812\,876&1341&23&659&13&200&792\,408& \multicolumn{2}{c}{?} & ?\\
\hline
\end{tabular}
\end{table}
\clearpage

The optimal Brinkhuis triples are not necessarily unique, and the list
also contains a case, $n=29$, where there exist optimal Brinkhuis
triples of both types. In general, several choices exist, which,
however, cannot be combined into an even larger triple. A list of
optimal $n$-Brinkhuis triples which at the same time are
minimal-length Brinkhuis $k_{\rm opt}$ triples is given below.

\begin{prop}
  The following sets of words generate optimal and minimal-length
  Brinkhuis triples:
\begin{itemize}
\item $n=13$, $k_{\rm p}=1$, $k_{\rm n}=0$, $k=1$:
\begin{equation}
\mathcal{G}_{13} = \{0121021201210\}
\end{equation}
\item $n=18$, $k_{\rm p}=0$, $k_{\rm n}=1$, $k=2$:
\begin{equation}
\mathcal{G}_{18} = \{012021020102120210\}
\end{equation}
\item $n=23$, $k_{\rm p}=1$, $k_{\rm n}=1$, $k=3$:
\begin{align}
\mathcal{G}_{23} = \{
&01210212021012021201210,\notag\\[-1ex]
&01210201021012021201210\}
\end{align}
\item $n=25$, $k_{\rm p}=1$, $k_{\rm n}=2$, $k=5$:
\begin{align}
\mathcal{G}_{25} = \{
&0121021202102012021201210,\notag\\[-1ex]
&0121020102101201021201210,\notag\\[-1ex]
&0121021201021012021201210\}
\end{align}
\item $n=29$, $k_{\rm p}=2$, $k_{\rm n}=2$, $k=6$:
\begin{align}
\mathcal{G}^{(1)}_{29} = \{
&01202120102012021020102120210,\notag\\[-1ex]
&01202120121012021012102120210,\notag\\[-1ex]
&01202102012101202120102120210,\notag\\[-1ex]
&01202120102012021012102120210\}
\end{align}
\item $n=29$, $k_{\rm p}=0$, $k_{\rm n}=3$, $k=6$:
\begin{align}
\mathcal{G}^{(2)}_{29} = \{
&01210201021201020121021201210,\notag\\[-1ex]
&01210201021202101201021201210,\notag\\[-1ex]
&01210201021202102012021201210\}
\end{align}
\item $n=30$, $k_{\rm p}=0$, $k_{\rm n}=4$, $k=8$:
\begin{align}
\mathcal{G}_{30} = \{
&012102010210120102012021201210,\notag\\[-1ex] 
&012102010212012102012021201210,\notag\\[-1ex]
&012102010212021020121021201210,\notag\\[-1ex] 
&012102120210120102012021201210\}
\end{align}
\item $n=33$, $k_{\rm p}=0$, $k_{\rm n}=6$, $k=12$:
\begin{align}
\mathcal{G}_{33} = \{
&012021020121012010212012102120210,\notag\\[-1ex] 
&012021020121021201021012102120210,\notag\\[-1ex] 
&012021020121021201210120102120210,\notag\\[-1ex] 
&012021201020120210121020102120210,\notag\\[-1ex] 
&012021201020121012021012102120210,\notag\\[-1ex] 
&012021201021012010212012102120210\}
\end{align}
\item $n=35$, $k_{\rm p}=2$, $k_{\rm n}=8$, $k=18$:
\begin{align}
\mathcal{G}_{35} = \{
&01202120102012102120121020102120210,\notag\\[-1ex] 
&01202120102101210201210120102120210,\notag\\[-1ex] 
&01202102010212010201202120102120210,\notag\\[-1ex] 
&01202102010212010201210120102120210,\notag\\[-1ex] 
&01202102012101201020121020102120210,\notag\\[-1ex] 
&01202102012101202120121020102120210,\notag\\[-1ex] 
&01202102012102120210121020102120210,\notag\\[-1ex] 
&01202120102012101201021012102120210,\notag\\[-1ex] 
&01202120102012102010210120102120210,\notag\\[-1ex] 
&01202120102120210201021012102120210\}
\end{align}
\item $n=36$, $k_{\rm p}=0$, $k_{\rm n}=16$, $k=32$:
\begin{align}
\mathcal{G}_{36} = \{
&012102010210120212010210121021201210,\notag\\[-1ex]
&012102010210120212012101201021201210,\notag\\[-1ex]
&012102010210121021201020121021201210,\notag\\[-1ex]
&012102010210121021201021012021201210,\notag\\[-1ex]
&012102010210121021202101201021201210,\notag\\[-1ex]
&012102010212012101201020121021201210,\notag\\[-1ex]
&012102010212012101201021012021201210,\notag\\[-1ex]
&012102010212012102120210121021201210,\notag\\[-1ex]
&012102010212021012010210121021201210,\notag\\[-1ex]
&012102010212021020102101201021201210,\notag\\[-1ex]
&012102120102101202120102012021201210,\notag\\[-1ex]
&012102120102101210201021012021201210,\notag\\[-1ex]
&012102120102101210212021012021201210,\notag\\[-1ex]
&012102120121012010212021012021201210,\notag\\[-1ex]
&012102120121012021201021012021201210,\notag\\[-1ex]
&012102120121020102120102012021201210\}
\end{align}
\item $n=40$, $k_{\rm p}=0$, $k_{\rm n}=24$, $k=48$:
\begin{align}
\mathcal{G}_{40} = \{
&0120210201210120102120210121020102120210,\notag\\[-1ex]
&0120210201210120212010210121020102120210,\notag\\[-1ex]
&0120210201210212010210120212012102120210,\notag\\[-1ex]
&0120210201210212012101201021012102120210,\notag\\[-1ex]
&0120210201210212012101202120121020120210,\notag\\[-1ex]
&0120210201210212012102010210120102120210,\notag\\[-1ex]
&0120210201210212012102010212012102120210,\notag\\[-1ex]
&0120210201210212012102012021012102120210,\notag\\[-1ex]
&0120210201210212012102012021020102120210,\notag\\[-1ex]
&0120210201210212021020120212012102120210,\notag\\[-1ex]
&0120212010201202101210201021012102120210,\notag\\[-1ex]
&0120212010201210120102012021012102120210,\notag\\[-1ex]
&0120212010201210120102101202120102120210,\notag\\[-1ex]
&0120212010201210120102120121020102120210,\notag\\[-1ex]
&0120212010201210120210201021012102120210,\notag\\[-1ex]
&0120212010201210120210201202120102120210,\notag\\[-1ex]
&0120212010201210120212010210120102120210,\notag\\[-1ex]
&0120212010201210212010201202120102120210,\notag\\[-1ex]
&0120212010201210212012101202120102120210,\notag\\[-1ex]
&0120212010210120102012101202120102120210,\notag\\[-1ex]
&0120212010212021012021201021012102120210,\notag\\[-1ex]
&0120212010212021012102010212012102120210,\notag\\[-1ex]
&0120212010212021012102120102012102120210,\notag\\[-1ex]
&0120212010212021020102120102012102120210\}
\end{align}
\item $n=41$, $k_{\rm p}=3$, $k_{\rm n}=31$, $k=65$:
\begin{align}
\mathcal{G}_{41} = \{
&01202102012102120210201202120121020120210,\notag\\[-1ex]
&01202120121012010201210201021012102120210,\notag\\[-1ex]
&01202120121021201021012010212012102120210,\notag\\[-1ex]
&01202102012101201021202101202120102120210,\notag\\[-1ex]
&01202102012101202120102012021012102120210,\notag\\[-1ex]
&01202102012101202120102012021020102120210,\notag\\[-1ex]
&01202102012102120102012101202120102120210,\notag\\[-1ex]
&01202102012102120121012010212012102120210,\notag\\[-1ex]
&01202102012102120121020120212012102120210,\notag\\[-1ex]
&01202120102012021012010210121020102120210,\notag\\[-1ex]
&01202120102012021012102010210120102120210,\notag\\[-1ex]
&01202120102012021012102010212012102120210,\notag\\[-1ex]
&01202120102012021012102012021020102120210,\notag\\[-1ex]
&01202120102012021012102120102012102120210,\notag\\[-1ex]
&01202120102012021012102120121020102120210,\notag\\[-1ex]
&01202120102012021020102120102012102120210,\notag\\[-1ex]
&01202120102012021020102120121020102120210,\notag\\[-1ex]
&01202120102012101201020120212012102120210,\notag\\[-1ex]
&01202120102012101202102010210120102120210,\notag\\[-1ex]
&01202120102012101202102010212012102120210,\notag\\[-1ex]
&01202120102012101202102012021012102120210,\notag\\[-1ex]
&01202120102012102120102012021012102120210,\notag\\[-1ex]
&01202120102012102120102101202120102120210,\notag\\[-1ex]
&01202120102012102120121012021012102120210,\notag\\[-1ex]
&01202120102012102120210201021012102120210,\notag\\[-1ex]
&01202120102012102120210201202120102120210,\notag\\[-1ex]
&01202120102101201020121012021012102120210,\notag\\[-1ex]
&01202120102101201021202101202120102120210,\notag\\[-1ex]
&01202120102120210121021201021012102120210,\notag\\[-1ex]
&01202120102120210201202120102012102120210,\notag\\[-1ex]
&01202120121012010201202120102012102120210,\notag\\[-1ex]
&01202120121012021012102010212012102120210,\notag\\[-1ex]
&01202120121012021012102120102012102120210,\notag\\[-1ex]
&01202120121012021020102120102012102120210\}
\end{align}
\end{itemize}
\end{prop}
\begin{proof}
  The proof that these are indeed special Brinkhuis triples consist of
  checking the conditions of definition~\ref{gentrip} explicitly. This
  has to be done by computer, as the number of symmetry-inequivalent
  composed words of length $3n$ that have to be checked for
  square-freeness is $k(2k^2+k_{\rm p})$, which gives $549\,445$ words
  of length $123$ for $\mathcal{G}_{41}$.  A Mathematica \cite{Wolf}
  program {\tt brinkhuistriples.m} that performs these checks
  accompanies this paper. This check is independent of the
  construction algorithm used to find the optimal triples.  In order
  to show that these triples are indeed optimal, one has to go through
  the algorithm outlined above. This has been done, giving the results
  of tables~\ref{trip1} and \ref{trip2}.
\end{proof}

The triple for $n=18$ is equivalent to the triple (\ref{EZtrip}) of
\cite{EZ}. The corresponding lower bounds on s are 
\begin{alignat}{2}
n=13, k=1:  \qquad& s\ge&  1^{1/12}&=1\notag\\[-1ex]
n=18, k=2:  \qquad& s\ge&  2^{1/17}&>1.041616\notag\\[-1ex]
n=23, k=3:  \qquad& s\ge&  3^{1/22}&>1.051204\notag\\[-1ex]
n=25, k=5:  \qquad& s\ge&  5^{1/24}&>1.069359\notag\\[-1ex]
n=29, k=6:  \qquad& s\ge&  6^{1/28}&>1.066083\notag\\[-1ex]
n=30, k=8:  \qquad& s\ge&  8^{1/29}&>1.074338\notag\\[-1ex]
n=33, k=12: \qquad& s\ge& 12^{1/32}&>1.080747\notag\\[-1ex]
n=35, k=18: \qquad& s\ge& 18^{1/34}&>1.088728\notag\\[-1ex]
n=36, k=32: \qquad& s\ge& 32^{1/35}&>1.104089\notag\\[-1ex]
n=40, k=48: \qquad& s\ge& 48^{1/39}&>1.104355\notag\\[-1ex]
n=41, k=65: \qquad& s\ge& 65^{1/40}&>1.109999
\end{alignat}
Apparently, the largest value of $n$ considered here yields the best
lower bound. This suggests that the bound can be systematically
improved by considering special Brinkhuis triples for longer words.

What about the restriction to spectial Brinkhuis triples? In general,
it is not clear what the answer is, but for the Brinkhuis triple pair
of \cite{EZ} it can easily be checked by computer that one cannot find
a shorter triple by lifting these restriction. In fact, this follows
from the following stronger result which is easier to check.

\begin{lemma}
  An $n$-Brinkhuis $(2,1,1)$-triple requires $n>17$.
\end{lemma}
\begin{proof}
  This can be checked by computer. The number of square-free words of
  length $n=17$ is $1044$. However, we do not need to check all
  $1044^4$ possibilities. Without loss of generality, we may restrict
  one of the four words to start with the letters $01$, leaving only
  $1044/6=174$ choices for this word. Furthermore, the two words in
  $\mathcal{B}^{(0)}$ may be interchanged, as well as the other two
  words; so it is sufficient to consider one order of words in both
  cases. No $n$-Brinkhuis $(2,1,1)$-triple was found for $n\le 17$.
\end{proof}

\section{Concluding remarks}

By enumerating square-free ternary words up to length $110$ and by
constructing generalised Brinkhuis triples, we improved both upper and
lower bounds for the number of ternary square-free words. The
resulting bounds for the exponential growth rate $s$ (\ref{s}) are
\begin{equation}
1.109999 < 65^{1/40} \le s\le 
{8\,416\,550\,317\,984}^{\frac{1}{108}}<1.317278.
\end{equation}
The main difficulty in improving the lower bound further is caused by
the combinatorial step in the algorithm to find optimal special
Brinkhuis triples. The data in tables~\ref{trip1} and \ref{trip2}
suggest that generators from the set $\mathcal{A}_{1}(n)$ (\ref{A1})
are more likely to provide optimal $n$-Brinkhuis triples for large $n$
than generators from the set $\mathcal{A}_{2}(n)$ (\ref{A2}). It would
be interesting to know whether, in principle, the lower bound obtained
in this way eventually converges to the actual value of $s$.

\section*{Acknowledgment}
The author would like to thank Jean-Paul Allouche for useful
discussions during a workshop at Oberwolfach in May 2001. The author
gratefully acknowledges comments from Shalosh B.~Ekhad and
Doron Zeilberger, who pointed out an error in a previous attempt to
improve the lower bound.

\end{document}